\documentclass[12pt]{iopart}

\usepackage{iopams}

\usepackage{graphicx}
\usepackage{amssymb}
\usepackage{amsthm}
\usepackage{latexsym}
\usepackage[colorinlistoftodos, textsize=tiny]{todonotes}
\usepackage{verbatim}
\usepackage{marginnote}
\usepackage[pagebackref]{hyperref} 		
\hypersetup{
	colorlinks,
	linkcolor=blue, 
	citecolor=blue, 
	linktocpage}

\theoremstyle{plain}
\newtheorem{theorem}{Theorem}[section]

\newtheorem{lemma}[theorem]{Lemma}

\theoremstyle{definition}
\newtheorem{definition}[theorem]{Definition}
\newtheorem{remark}[theorem]{Remark}
\newtheorem{assumption}[theorem]{Assumption}

\def\clo#1{\overline{#1}}
\def\text#1{\mbox{#1}}

\newcommand{\ra}{\rangle}
\newcommand{\la}{\langle}

\newcommand{\Sop}{\mathcal{S}}
\newcommand{\Cop}{\mathcal{C}}

\newcommand{\Aop}{\mathcal{A}}
\newcommand{\Uop}{\mathcal{U}}

\newcommand{\g}{\mathrm{g}}
\newcommand{\grad}{\text{grad}}
\newcommand{\dVol}{\text{dVol}}
\newcommand{\CD}{\text{CD}}
\newcommand{\dS}{\text{dS}}

\usepackage{color}
\definecolor{hotcolor}{rgb}{1,0,0}

\usepackage[letterpaper,margin=1.15in]{geometry}

\begin{document}

\title{Multiwave imaging in an enclosure with variable wave speed}

\author{Sebasti\'{a}n Acosta$^{1}$ and Carlos Montalto$^2$}
\address{$^1$ Department of Pediatrics -- Cardiology, Baylor College of Medicine, TX, USA}

\address{$^2$ Department of Mathematics, University of Washington, Seattle, WA, USA}

\eads{\mailto{sacosta@bcm.edu} and \mailto{montcruz@uw.edu}}

\begin{abstract}
In this paper we consider the mathematical model of thermo- and photo-acoustic tomography for the recovery of the initial condition of a wave field from knowledge of its boundary values. Unlike the free-space setting, we consider the wave problem in a region enclosed by a surface where an impedance boundary condition is imposed. This condition models the presence of physical boundaries such as interfaces or acoustic mirrors which reflect some of the wave energy back into the enclosed domain. By recognizing that the inverse problem is equivalent to a statement of boundary observability, we use control operators to prove the unique and stable recovery of the initial wave profile from knowledge of boundary measurements. Since our proof is constructive, we explicitly derive a solvable equation for the unknown initial condition. This equation can be solved numerically using the conjugate gradient method. We also propose an alternative approach based on the stabilization of waves. This leads to an exponentially and uniformly convergent Neumann series reconstruction when the impedance coefficient is not identically zero. In both cases, if well-known geometrical conditions are satisfied, our approaches are naturally suited for variable wave speed and for measurements on a subset of the boundary.
\end{abstract}


\submitto{\IP}
% \maketitle

%\listoftodos

%%%%%%%%%%%%%%%%%%%%%%%%%%%%%%%%%%%%%%%%%%%%%%%%%%%%%%%%%%%%%%%%%
%%%%%%% NEW SECTION %%%%%%%%%%%%%%%%%%%%%%%%%%%%%%%%%%%%%%%%%%%%%
%%%%%%%%%%%%%%%%%%%%%%%%%%%%%%%%%%%%%%%%%%%%%%%%%%%%%%%%%%%%%%%%%
\section{Introduction} \label{Section:Intro}

\textit{Photoacoustic Tomography} (PAT) and \textit{Thermoacoustic Tomography} (TAT) are medical imaging modalities that combine the high contrast properties of electromagnetic waves, with the high resolution properties of acoustic waves. When a short electromagnetic pulse (laser/microwave) is used to irradiate a biological tissue, the photoacousic/thermoacoustic effect results in emission of acoustic signals that can be measured outside the object by wide-band ultrasonic transducers \cite{Wang-Wu-2007,Wang-2009,Beard-2011, Kuchment-2008, Anastasio-2011, Wang-2012}. The goal of PAT and TAT is to obtain an image of the electromagnetic absorption properties of the tissue. These properties are highly related with the molecular structure of the tissue, and hence could reveal pathological conditions \cite{Wang-2003, Zhang-2007, Shan-2008, Pramanik-Wang-2009}. The process of recovering the actual image of the tissue from the external data involves two steps. The first step consists of recovering the internal absorbed radiation from knowledge of the acoustic signal measured in the exterior. In both PAT and TAT, the first step is the same and is modeled as the reconstruction of an initial condition of a wave field from data measured on the boundary of the region. \cite{Beard-2011, Bal-2011}. 

In most of the mathematical analysis of TAT and PAT, there is the assumption that the acoustic waves propagate in open space \cite{Agranovsy-2007, Hristova-et-al-2008, Ste-Uhl-2009-01, Stefanov-Uhlmann-2011, Homan-2013, Tittelfitz-2012, Anastasio-2007}. However, in practice this is not always a valid assumption \cite{Wang-Yang-2007}. In some cases, the acoustic waves interact with boundaries/interfaces such as the skull, detectors, or the interface with air \cite{Jin-2008, Schoonover-Anas-2011, Nie-Wang-2012, Schoonover-2012,Huang-Nie-2012}. In other cases, certain boundaries are purposely introduced in the experimental setup to improve performance. For example, the bioengineering group at University College London observed that better imaging resolution and an artificially increased detection aperture can be attained by partially or fully enclosing the target region with acoustic mirrors \cite{Cox-Arridge-Beard-2007,Cox-Beard-2009,Ellwood2014}. Another example was proposed in \cite{Huang-2013}, where a $45^\circ$ angle acoustic reflector is used to improve the limited detection view of linear arrays. Improved image quality on  phantom and ex-vivo experiments were observed. Moreover, their method can be implemented by having the acoustic reflectors on the boundary of an immersion tank. In such cases, the detectors will form an half-enclosed space for the object, enhancing angle detection and potentially improving the image quality substantially.  In the mathematical literature this problem was first analyzed in \cite{kunyansky-holman-cox-2013,holman-kunyansky-2014}, where reconstruction algorithms are obtained for an isotropic medium in rectangular cavities. In \cite{kunyansky-holman-cox-2013} a method based on the fast Fourier transform was presented to study the case of constant speed of a rectangular domain. Holman and Kunyansky \cite{holman-kunyansky-2014} used an eigenfunction expansion method to obtain an approximate reconstruction that improves as the measurement time increases.

In our work, we recognize that this inverse problem is equivalent to a statement of observability in the theory of control for partial differential equations (PDEs). Under the so-called geometric control condition (see Assumption \ref{Assump.001}), uniqueness and stability follows from the work of Bardos, Lebeau and Rauch  \cite{Bar-Leb-Rau-1992}. This control approach is naturally suited for variable anisotropic wave speed and for measurements on a subset of the boundary. Also, to model the presence of a physical boundary, such as detectors or acoustic mirrors, we consider the problem in a region enclosed by a surface where an impedance boundary condition is imposed. The impedance coefficient is allowed to vary over all non-negative values. Zero values of the impedance coefficient model acoustically hard surfaces, intermediate values model absorption, and  relatively high values approximate acoustically soft boundaries (such as the interface with air).

We present two different reconstruction schemes for this problem. The first one is based on the proof of uniqueness and stability. We explicitly derive a solvable equation for the unknown initial condition. This equation can be solved using the conjugate gradient iterative method. The second approach is based on a Neumann series expansion introduced by Stefanov and Uhlmann \cite{stefanov2009thermoacoustic} for the case of free space propagation. We show that the Neumann series reconstruction is valid when the geometric control condition holds for the region of the boundary where the impedance coefficient is positive.  It is worth mentioning that this work was developed simultaneously to the work of Stefanov and Yang \cite{stefanov-yang-2014} where they present an average time reversal approach to obtain a Neumann series reconstruction in the case when the impedance coefficient is identically zero.

The structure of this paper is as follows. Section \ref{Section:Formulation} contains the mathematical formulation of the forward acoustic problem where we incorporate the impedance boundary condition and the possibility of anisotropic variable media modeled in geometric terms. Anisotropy is commonly encountered in biological media such as bones and muscles, and the Riemannian geometric formulation is a mathematically convenient manner to model this type of media \cite{GLLT-2004,stefanov2009thermoacoustic}. In Section \ref{Section:Control} we briefly review some tools of control theory which we employ in the remainder of the paper. These tools are naturally suited for variable wave speed and for measurements on a subset of the boundary. Section \ref{Section:UniqStab} contains the main abstract results concerning uniqueness and stability of the inversion, whereas Sections \ref{Section:ReconstCGM} and \ref{Section:ReconstNeumann} describe our proposed constructive algorithms to recover the initial acoustic profile from the boundary measurement. Finally, we provide a brief discussion of these results and of future work in Section \ref{Section:Discussion}.

%%%%%%%%%%%%%%%%%%%%%%%%%%%%%%%%%%%%%%%%%%%%%%%%%%%%%%%%%%%%%%%%%
%%%%%%% NEW SECTION %%%%%%%%%%%%%%%%%%%%%%%%%%%%%%%%%%%%%%%%%%%%%
%%%%%%%%%%%%%%%%%%%%%%%%%%%%%%%%%%%%%%%%%%%%%%%%%%%%%%%%%%%%%%%%%
\section{Mathematical Formulation} \label{Section:Formulation}

In this paper we study the thermoacoustic tomography problem in the presence of heterogeneous media and an enclosing surface. More precisely we assume anisotropy being modeled in geometric terms. For dimension $n \geq 2$, let $\Omega \subset \mathbb{R}^{n}$ be a smooth, bounded and simply connected domain with boundary $\partial \Omega$. Let $\g$ be a Riemannian metric defined on $\clo{\Omega}$. The propagation of acoustic waves is governed by the following system,
\numparts 
\begin{eqnarray}
\partial_{t}^{2} {u} - \Aop u = 0 \qquad && \mbox{in} \quad (0,\tau) \times \Omega \label{Eqn:001f}  \\
u = u_{0} \quad \mbox{and} \quad \partial_{t} u = 0 \qquad && \mbox{on} \quad \{ t = 0 \} \times \Omega \label{Eqn:002f} \\
\partial_{\nu} u + \lambda \partial_{t} u  = 0 \qquad && \mbox{on} \quad (0,\tau) \times \partial \Omega \label{Eqn:003f}
\end{eqnarray}
\endnumparts
where $u_{0} \in H^{1}(\Omega)$ is the unknown initial acoustic profile. Here $\Aop$ denotes a second-order elliptic differential operator of the form
\begin{eqnarray} 
\Aop u = c^2 \Delta_{\g} u - q u, \label{Eqn:DiffOp}
\end{eqnarray}
where $c: \Omega \to \mathbb{R}$ is the wave speed, $\Delta_{\g}$ is the Laplace-Beltrami operator associated with the Riemannian metric $\g$ and $q: \Omega \to \mathbb{R}$ is a potential.  Also,  $\lambda :\partial \Omega \to \mathbb{R}$ denotes the impedance coefficient at the boundary $\partial \Omega$. We assume that $\lambda$ is a non-negative function. Physically, $\lambda = 0$ models an acoustically hard surface and $\lambda \to \infty$ approximates an acoustically soft boundary. In general we allow $\lambda$ to vary on the boundary $\partial \Omega$ to model the presence or absence of detectors, or the inhomogeneous nature of an enclosing surface. The length of the observation window of time is given by $\tau < \infty$ which is defined in the next section. In (\ref{Eqn:003f}), the symbol $\partial_{\nu}$ denotes the conormal derivative with respect to the metric $\g$ at the boundary $\partial \Omega$. All the coefficients in (\ref{Eqn:001f})-(\ref{Eqn:003f}) are assumed to be sufficiently smooth.

In order to consider partial measurements, we divide the boundary as the disjoint union $\partial \Omega = \Gamma \cup (\partial \Omega \setminus \Gamma)$ and we assume that $\Lambda u_0 = u|_{(0,\tau) \times \Gamma}$ constitutes our measurements. Then, the goal of the thermoacoustic tomography problem is to find the initial profile $u_{0}$ from knowledge of $\Lambda u_0$. This is a challenging problem with intricate dependencies on the nature of the manifold $(\Omega,c^{-2} \g)$, the partial boundary $\Gamma$ and the time interval $(0,\tau)$. We recognize that, in mathematical terms, the thermoacoustic tomography problem coincides with a problem of \textit{boundary observability}, which is one of the central concepts of control theory for PDEs \cite{Glow-Lions-He-2008,Bar-Leb-Rau-1992,GLLT-2004}. Hence, the objective of this paper is to constructively employ the tools of control theory to solve the thermoacoustic tomography problem in an enclosed domain with heterogeneous sound speed.

%%%%%%%%%%%%%%%%%%%%%%%%%%%%%%%%%%%%%%%%%%%%%%%%%%%%%%%%%%%%%%%%%
%%%%%%% NEW SECTION %%%%%%%%%%%%%%%%%%%%%%%%%%%%%%%%%%%%%%%%%%%%%
%%%%%%%%%%%%%%%%%%%%%%%%%%%%%%%%%%%%%%%%%%%%%%%%%%%%%%%%%%%%%%%%%
\section{Background on Control Theory} \label{Section:Control}

Our approach relies on exact boundary controllability for the wave equation. Hence, the purpose of this section is to review some background and define the notation concerning control theory. For details, we refer the reader to \cite{Glow-Lions-He-2008,Bar-Leb-Rau-1992,GLLT-2004,Lions-Review-1988}. The first point is to ensure that the manifold $(\Omega,c^{-2} \g)$, on which the wave equation is defined, yields exact boundary controllability. 

Following Bardos, Lebeau and Rauch \cite{Bar-Leb-Rau-1992}, we assume that our problem enjoys the \textit{geometric control condition} (GCC) for the Riemannian manifold $(\Omega, c^{-2} \g)$ with only a portion $\Gamma$ of the boundary $\partial \Omega$ being accessible for control and observation. We assume that $\Gamma$ is a smooth and simply connected domain relative to $\partial \Omega$. In this paper, we assume that $\Gamma$ satisfies the following.

\begin{assumption}[\textbf{GCC}] \label{Assump.001} 
The geodesic flow of $(\Omega, c^{-2} \g)$ reaches $\Gamma$ after possible reflections on $\partial \Omega \setminus \Gamma$ in finite time $\tau$, and $\{\lambda > 0 \} \subset \Gamma$. In other words, there exists $\tau < \infty$ such that any geodesic ray, originating from any point in $\Omega$ at $t=0$, eventually reaches $\Gamma$ in a non-diffractive manner (after possible geometrical reflections on $\partial \Omega \setminus \Gamma$) before time $t=\tau$. We refer to $\tau$ as the exact controllability time.
\end{assumption}

Throughout the paper $\{ \lambda > 0 \}$ denotes the set $\{ x \in \partial \Omega: \lambda(x) > 0\}$, similarly $\{\lambda = 0\}$ denotes the set $\{x \in \partial \Omega: \lambda(x) = 0\}$. The condition $\{\lambda>0\} \subset\Gamma$ is a natural exact controllability condition for the initial problem (\ref{Eqn:001c})-(\ref{Eqn:004c}) with Robin boundary conditions. Concerning the inverse problem, this condition guarantees uniqueness, stability and reconstruction using a conjugate gradient method (see Sections \ref{Section:UniqStab}-\ref{Section:ReconstCGM}). However for a Neumann series reconstruction we need to assume that $ \{\lambda>0\} = \Gamma$ (see Section \ref{Section:ReconstNeumann}).

We work with the standard Sobolev spaces based on square-integrable functions. The associated inner-product extends naturally as the duality pairing between functionals and functions. We should interpret the Hilbert space $H^{0}(\Omega)$ with the inner-product appropriately weighted by $c^{-2}|\det\g|^{1/2}$ so that $c^2\Delta_{\g}$ is formally self-adjoint with respect to the duality pairing of $H^{0}(\Omega)$. Sometimes this is written as $H^{0}(\Omega; c^{-2}\text{dVol})$, where $\text{dVol}(x) = |\det \g|^{1/2}\text{d}x$, however to simplify notation and since there is no risk of confusion we will only write it as $H^{0}(\Omega)$. The same weight, $c^{-2}|\det\g|^{1/2}$,  is incorporated in the inner product for $H^{0}(\partial \Omega)$. This implies that the differential operator $\Aop$ is formally self-adjoint.

Now we consider the following auxiliary problem. Given $\zeta \in H^{0}((0,\tau) \times \Gamma)$, find the generalized solution $\xi \in C^{k}([0,\tau];H^{1-k}(\Omega))$ for $k=0,1$ of the following problem
\numparts
\begin{eqnarray} 
\partial_{t}^{2} \xi - \mathcal{A} \xi = 0 \qquad && \mbox{in} \quad (0,\tau) \times \Omega, \label{Eqn:001c}  \\
\xi = 0 \quad \mbox{and} \quad \partial_{t} \xi = 0 \qquad && \mbox{on} \quad \{ t = 0 \} \times \Omega, \label{Eqn:002c} \\
\partial_{\nu} \xi + \lambda \partial_{t} \xi    = \zeta \qquad && \mbox{on} \quad (0,\tau) \times \Gamma, \label{Eqn:003c} \\
\partial_{\nu} \xi  = 0 \qquad && \mbox{on} \quad (0,\tau) \times (\partial \Omega \setminus \Gamma). \label{Eqn:004c} 
\end{eqnarray}
\endnumparts
Since we assume that $\lambda \geq 0$ in the boundary conditions (\ref{Eqn:003c}) then the above is a well-posed problem satisfying a stability estimate of the following form,
\begin{eqnarray*}
\| \xi \|_{C^{k}([0,\tau];H^{1-k}(\Omega))} + \| \text{grad}_{\g} \, \xi \|_{C^{k}([0,\tau];H^{-k}(\Omega))} \lesssim \| \zeta \|_{H^{0}((0,\tau) \times \Gamma)},
\end{eqnarray*}
for $k=0,1$. For details see \cite{Glow-Lions-He-2008,Bar-Leb-Rau-1992} and references therein. Here and in the rest of the paper, the symbol $\lesssim$ means inequality up to a positive constant. This gives rise to the following definition of a bounded operator.

\begin{definition}[\textbf{Solution Operator}] \label{Def:SolOp}
Let the \textit{solution} operator
\begin{eqnarray*}
\Sop : H^{0}((0,\tau) \times \Gamma) \to H^{0}(\Omega), 
\end{eqnarray*}
be given by the map $\zeta \mapsto \partial_{t} \xi |_{t=\tau}$ where $\xi$ is the solution of (\ref{Eqn:001c})-(\ref{Eqn:004c}).
\end{definition}

The goal of the control problem is to prove that the operator $\Sop$ is surjective. In other words, given arbitrary $\phi \in H^{0}(\Omega)$, the goal is to find a boundary condition $\zeta \in H^{0}((0,\tau) \times \Gamma)$ to drive the solution $\xi$ of (\ref{Eqn:001c})-(\ref{Eqn:004c}) from vanishing Cauchy data at $t=0$ to the desired Cauchy data $( \xi, \partial_{t} \xi) = (0, \phi)$ at time $t=\tau$. For our own reference hereafter, we state the well-posedness of this control problem as a theorem. A proof obtained under the \textit{geometric control condition} is found in \cite{Bar-Leb-Rau-1992}. See also the extension in \cite{Burq-1997} to less regular domains and coefficients.

\begin{theorem}[\textbf{Exact Controllability}] \label{Thm.Control}
Let the geometric control condition \ref{Assump.001} hold. Then for any function $\phi \in H^{0}(\Omega)$, there is a boundary control $\zeta \in H^{0}((0,\tau) \times \Gamma)$ so that the solution $\xi$ of (\ref{Eqn:001c})-(\ref{Eqn:004c}) satisfies
\begin{eqnarray*}
(\xi , \partial_{t} \xi) = (0, \phi) \qquad \text{at time $t = \tau$.}
\end{eqnarray*}

Among all such boundary controls there exists $\zeta_{\rm min}$ which is uniquely determined by $\phi$ as the minimum norm control and satisfies the following stability condition
\begin{eqnarray*}
& \| \zeta_{\rm min} \|_{H^{0}((0,\tau) \times \Gamma)} \leq C \| \phi \|_{H^{0}(\Omega)}
\end{eqnarray*}
for some positive constant $C = C(\Omega,\Gamma,\Aop,\tau)$.
\end{theorem}

Along with the above theorem, we readily obtain the following definition of the control operator that we will use in the rest of the paper.

\begin{definition}[\textbf{Control Operator}] \label{Def:CSOpe}
Let the \textit{control} operator
\begin{eqnarray*}
& \Cop : H^{0}(\Omega) \to H^{0}((0,\tau) \times \Gamma), 
\end{eqnarray*}
be given by the map $\phi \mapsto \zeta_{\rm min}$ where $\zeta_{\rm min}$ is defined in Theorem \ref{Thm.Control}.
\end{definition}

This leads to the following lemma which states the basic properties for the \textit{control} operator. The proof is a direct consequence of Theorem \ref{Thm.Control}.

\begin{lemma} \label{Lemma.007}
The \textit{control} operator given by Definition \ref{Def:CSOpe} is bounded. Moreover, we have that $\phi = \Sop \Cop \phi$, and the following estimates,
\begin{eqnarray}
\| \phi \|_{H^{0}(\Omega)} \lesssim \| \Cop \phi \|_{H^{0}((0,\tau) \times \Gamma)}  \lesssim \| \phi \|_{H^{0}(\Omega)},  \label{Eqn:Stab02}
\end{eqnarray}
hold for all $\phi \in H^{0}(\Omega)$.
\end{lemma}

%%%%%%%%%%%%%%%%%%%%%%%%%%%%%%%%%%%%%%%%%%%%%%%%%%%%%%%%%%%%%%%%%
%%%%%%% NEW SECTION %%%%%%%%%%%%%%%%%%%%%%%%%%%%%%%%%%%%%%%%%%%%%
%%%%%%%%%%%%%%%%%%%%%%%%%%%%%%%%%%%%%%%%%%%%%%%%%%%%%%%%%%%%%%%%%
\section{Uniqueness and Stability using Control Theory} \label{Section:UniqStab}

In this Section we state the inverse problem for thermoacoustic tomography and we also introduce and prove our first main result. We start by recalling the forward problem (\ref{Eqn:001f})-(\ref{Eqn:003f}) which is well-posed with a solution $u \in C^{k}([0,\tau];H^{1-k}(\Omega))$ for $k=0,1$ for an initial condition $u_{0} \in H^{1}(\Omega)$.

\begin{definition}[\textbf{Inverse Problem}] \label{Def.InvProblem}
Under the geometric control condition \ref{Assump.001}, find the unknown initial condition $u_{0}$ of the forward problem (\ref{Eqn:001f})-(\ref{Eqn:003f}) from knowledge of the partial Dirichlet data $\Lambda u_0 = u|_{(0,\tau) \times \Gamma}$.
\end{definition}

Now we are ready to reduce the inverse problem to a solvable equation using duality arguments. First, let us introduce an auxiliary \textit{time-reversal} operator $\Uop : H^{0}((0,\tau) \times \Omega) \to H^{0}((0,\tau) \times \Omega)$ given by
\begin{eqnarray}
(\Uop v)(t,x) = v(\tau - t,x).  \label{Eqn.TimeReversalOp}
\end{eqnarray}
The operator $\Uop$ is clearly unitary and self-adjoint. The same is true for this operator defined as $\Uop : H^{0}((0,\tau) \times \Gamma) \to H^{0}((0,\tau) \times \Gamma)$. Also notice that $\Uop \partial_{t} v = - \partial_{t} \Uop v$ for functions that posses a weak derivative in time.

In what follows, we will evaluate the duality pairing between the terms in equation (\ref{Eqn:001f})-(\ref{Eqn:003f}) against $\psi = \Uop \xi$ where $\xi$ is the solution of (\ref{Eqn:001c})-(\ref{Eqn:004c}) with $\zeta = \Cop \phi$ and $\phi \in H^{0}(\Omega)$. Now, since $u$ is a weak solution of (\ref{Eqn:001f})-(\ref{Eqn:003f}) then it satisfies,
\begin{eqnarray*} 
\la \partial_{t} u , \psi \ra_{\Omega} \big|_{t=0}^{t=\tau} & = \la \partial_{\nu} u , \psi \ra_{(0,\tau) \times \partial \Omega} - \la \text{grad}_{\g} u , \text{grad} \, \psi \ra_{(0,\tau) \times \Omega} \\ 
& \qquad - \la c^{-2} q u , \psi \ra_{(0,\tau) \times \Omega} + 
\la \partial_{t} u , \partial_{t} \psi \ra_{(0,\tau) \times \Omega}.
\end{eqnarray*}
Similarly, $\psi = \Uop \xi$ is also a weak solution to the wave equation (\ref{Eqn:001c}), then it satisfies,
\begin{eqnarray*}
\la  u , \partial_{t} \psi \ra_{\Omega} \big|_{t=0}^{t=\tau} & = \la  u , \partial_{\nu} \psi \ra_{(0,\tau) \times \partial \Omega} - \la \text{grad}_{\g} u , \text{grad} \, \psi \ra_{(0,\tau) \times \Omega} \\
& \qquad - \la u , c^{-2} q \psi \ra_{(0,\tau) \times \Omega}  + 
\la \partial_{t} u , \partial_{t} \psi \ra_{(0,\tau) \times \Omega}.
\end{eqnarray*}
Subtracting the above identities and using the Cauchy data $u = u_{0}$ and $\partial_{t}u = 0$ at $t = 0$, and the vanishing Cauchy data $\psi = \partial_{t} \psi = 0$ at time $t = \tau$, and the fact that $\partial_{t} \psi |_{t=0} = - \phi $,  we obtain,
\begin{eqnarray*}
 \la u_{0} , \phi \ra_{\Omega} & = \la  u , \partial_{\nu} \psi \ra_{(0,\tau) \times \partial \Omega} - \la \partial_{\nu} u , \psi \ra_{(0,\tau) \times \partial \Omega}  \\
& = \la  u , \partial_{\nu} \psi - \lambda \partial_{t} \psi \ra_{(0,\tau) \times \partial \Omega} - \la \partial_{\nu} u + \lambda \partial_{t}u , \psi \ra_{(0,\tau) \times \partial \Omega}.
\end{eqnarray*}
Now notice that $\partial_{\nu} \psi = \partial_{\nu} \Uop \xi = \Uop \partial_{\nu} \xi$ and $\lambda \partial_{t} \psi = \lambda \partial_{t} \Uop \xi = - \Uop \lambda \partial_{t} \xi$ on $(0,\tau) \times \partial \Omega$. Therefore, from the boundary conditions (\ref{Eqn:003f}) and (\ref{Eqn:003c})-(\ref{Eqn:004c}), we arrive at
\begin{eqnarray*}
\la u_{0} , \phi \ra_{\Omega} =  \la  \Lambda u_0 , \Uop \zeta \ra_{(0,\tau) \times \Gamma }  = \la  \Lambda u_0 , \Uop \Cop \phi \ra_{(0,\tau) \times \Gamma }, \quad \text{for all $\phi \in H^{0}(\Omega)$,}
\end{eqnarray*}
where $\Lambda u_0=u |_{(0,\tau) \times \Gamma}$ is the measurement at the partial boundary $\Gamma$.

As a result, the unknown initial condition $u_{0}$ is explicitly recovered as follows,
\begin{eqnarray}
u_{0} = \Cop^{*} \Uop \Lambda u_0, \label{Eqn.Solution}
\end{eqnarray}
where $\Cop^{*} : H^{0}((0,\tau) \times \Gamma) \to H^{0}(\Omega)$ is the adjoint of $\Cop : H^{0}(\Omega) \to H^{0}((0,\tau) \times \Gamma)$ which is a bounded operator. This proves the first of our main results.

\begin{theorem}[\textbf{Uniqueness and Stability}]
The inverse problem of thermoacoustic tomography, as stated in Definition \ref{Def.InvProblem}, is uniquely solvable. Moreover, the following stability  estimate
\begin{eqnarray*}
\| u_{0} \|_{H^{0}(\Omega)} \leq C \| \Lambda u_0 \|_{H^{0}((0,\tau) \times \Gamma)}, \label{Eqn.Stability}
\end{eqnarray*}
holds for a positive constant $C = C(\Omega, \Gamma, \Aop,\tau)$.
\end{theorem}

\begin{remark}
From the details in \cite{Bar-Leb-Rau-1992}, we have that the control operator $\Cop$ is bounded on Sobolev spaces of any scale of regularity (positive and negative). Therefore, we easily obtain from (\ref{Eqn.Solution}) a stability estimate of the form $\| u_{0} \|_{H^{1}(\Omega)} \lesssim \| \Lambda u_{0} \|_{H^{1}((0,\tau) \times \Gamma)}$ which matches the level of regularity assumed in the forward problem (\ref{Eqn:001f})-(\ref{Eqn:003f}).
\end{remark}

%%%%%%%%%%%%%%%%%%%%%%%%%%%%%%%%%%%%%%%%%%%%%%%%%%%%%%%%%%%%%%%%%
%%%%%%% NEW SECTION %%%%%%%%%%%%%%%%%%%%%%%%%%%%%%%%%%%%%%%%%%%%%
%%%%%%%%%%%%%%%%%%%%%%%%%%%%%%%%%%%%%%%%%%%%%%%%%%%%%%%%%%%%%%%%%
\section{Reconstruction using Conjugate Gradient Method} \label{Section:ReconstCGM}

In order to take full advantage of the recovery equation (\ref{Eqn.Solution}), we now proceed to express the adjoint-control operator $\Cop^{*}$ in terms of applicable operators. Hence, we remind the reader of the following facts from control theory. For a proof see \cite[Ch. 6]{Glow-Lions-He-2008} and \cite{Bar-Leb-Rau-1992,Lions-Review-1988}.

\begin{lemma} \label{Lemma.009}
Let the mapping $\Sop^{*} : H^{0}(\Omega) \to H^{0}((0,\tau) \times \Gamma)$ be given by $z \mapsto w|_{(0,\tau) \times \Gamma}$ where $w \in C^{k}([0,\tau] ; H^{-k}(\Omega))$ for $k=0,1$ is the generalized solution of the following \textit{time-reversed} problem,
\numparts
\begin{eqnarray} 
\partial_{t}^{2} w - \mathcal{A} w = 0 \qquad && \mbox{in} \quad (0,\tau) \times \Omega \label{Eqn:001a}  \\
w = z \quad \mbox{and} \quad \partial_{t} w = 0 \qquad && \mbox{on} \quad \{ t = \tau \} \times \Omega \label{Eqn:002a} \\
\partial_{\nu} w - \lambda \partial_{t} w    = 0 \qquad && \mbox{on} \quad (0,\tau) \times \partial \Omega. \label{Eqn:003a}
\end{eqnarray}
\endnumparts
Then, as the notation suggests, $\Sop^{*}$ is the adjoint of the solution operator $\Sop$ from Definition \ref{Def:SolOp}. Moreover, if the geometrical control condition \ref{Assump.001} holds, then we have:
\begin{itemize}
\item[1.] The operator $(\Sop \Sop^{*}) :H^{0}(\Omega) \to H^{0}(\Omega)$ is coercive.
\item[2.] The control operator $\Cop$ is given by $\Cop = \Sop^{*} (\Sop\Sop^{*})^{-1} $. 
\item[3.] An equation of the type $(\Sop \Sop^{*}) \phi = \beta$ can be solved using the conjugate gradient method.
\end{itemize}
\end{lemma}

Our second main result follow directly from Lemma \ref{Lemma.009} and the explicit formula for $u_0$ obtained in (\ref{Eqn.Solution}).

\begin{theorem}[\textbf{Reconstruction Method 1}] \label{Thm.Reconstruction}
The adjoint of the control operator $\Cop$ satisfies $\Cop^{*} = (\Sop \Sop^{*})^{-1} \Sop$. Hence, the unique solution $u_{0}$ of the thermoacoustic tomography problem is given by,
\begin{eqnarray}
u_{0} = (\Sop \Sop^{*})^{-1} \Sop \, \Uop \, \Lambda u_0, \label{Eqn.SolutionAgain}
\end{eqnarray}
where $\Uop$ is the time-reversal operator defined in (\ref{Eqn.TimeReversalOp}) and $(\Sop \Sop^{*})$ can be inverted using the conjugate gradient iterative method with convergence in $H^{0}(\Omega)$.
\end{theorem}

%%%%%%%%%%%%%%%%%%%%%%%%%%%%%%%%%%%%%%%%%%%%%%%%%%%%%%%%%%%%%%%%%
%%%%%%% NEW SECTION %%%%%%%%%%%%%%%%%%%%%%%%%%%%%%%%%%%%%%%%%%%%%
%%%%%%%%%%%%%%%%%%%%%%%%%%%%%%%%%%%%%%%%%%%%%%%%%%%%%%%%%%%%%%%%%
\section{Reconstruction by Neumann Series}
\label{Section:ReconstNeumann}

In this section we consider an alternative reconstruction using a Neumann series expansion in the case when $\Gamma = \{ \lambda > 0 \}$ is non empty and satisfies the GCC  \ref{Assump.001}. In other words, we show that if $\Gamma$ satisfies the GCC \ref{Assump.001} with $ \{ \lambda > 0 \} \subset \Gamma$ replaced by  $\Gamma = \{ \lambda > 0 \}$, then we can recover $u_0$ from a Neumann series applied to a backprojection of the data restricted to $\Gamma$. Our analysis follows mainly the stabilization of waves in the work by Bardos, Lebeau and Rauch \cite[Sect. 5]{Bar-Leb-Rau-1992} and uses in a non-trivial way the energy decay  guaranteed by $\{\lambda >0\} = \Gamma$. For the case $\lambda \equiv 0$ the construction presented here will not work, even in the case of full data, because a certain error operator $K$ will fail to be a contraction due to the fact that total energy in conserved. For this case, we refer the reader to  \cite{stefanov-yang-2014} were an average time reversal approach is used to avert this difficulty and obtain a Neumann series expansion when $\lambda \equiv 0$.

We use a similar backprojection operator as the one introduced by Stefanov and Uhlmann \cite{stefanov2009thermoacoustic} for the free-space setting. See also \cite{qian2011efficient}. However in our case, we can use the information of the Neumann derivative of the wave in $\{\lambda >0\}$ to make the error operator $K$ satisfy Neumann boundary condition in all of $\partial \Omega$. This avoids using a harmonic extension at $t = \tau$ and allows us to obtain a Neumann series reconstruction with information only on the partial boundary $\Gamma$.

We introduce some notation for energy spaces with Robin boundary condition and describe some properties of the evolution operator for the Cauchy problem. 

\subsection{Energy Spaces and Evolution Operators}
Consider the initial Cauchy value problem
\numparts 
\begin{eqnarray}
\partial_{t}^{2} {u} - \Aop u = 0 \qquad && \mbox{in} \quad (0,\tau) \times \Omega \label{Eqn:001f-Cauchy}  \\
u = u_{0} \quad \mbox{and} \quad \partial_{t} u = u_1 \qquad && \mbox{on} \quad \{ t = 0 \} \times \Omega \label{Eqn:002f-Cauchy} \\
\partial_{\nu} u + \lambda \partial_{t} u  = 0 \qquad && \mbox{on} \quad (0,\tau) \times \partial \Omega \label{Eqn:003f-Cauchy}
\end{eqnarray}
\endnumparts 
for $\lambda$ a non-negative smooth function. For this subsection we allow $\lambda \equiv 0$.

Define the energy of $u$ at time $t$ by
\begin{equation} \label{eq:Energy}
E(t,u) = \frac{1}{2}\int_{\Omega } |\grad_\g u|^2 + c^{-2}q|u|^2 + c^{-2}|u_t|^2 \dVol.
\end{equation}

Denote by $\CD_R^\infty(\Omega)$ the set of Cauchy data related to Robin boundary condition $(u_0,u_1)$ for which there exist a unique smooth solutions of the initial problem (\ref{Eqn:001f-Cauchy})-(\ref{Eqn:003f-Cauchy}) in $(0,\tau)\times \Omega$.
 
%In particular, $\CD_R^\infty(\Omega)$ contains all Cauchy data $(u_0,u_1)$ satisfying compatibility conditions of all orders, i.e., 
%\[
%\CD_R^\infty(\Omega)  \subset  \{(u_0,u_1) \in  C^\infty(\overline{\Omega})^2: \partial_t^k(\partial_{\nu} + \lambda \partial_t) u |_{\{t=0\} \times\partial\Omega} =0,  k=0,1, \ldots \}.
%\]

Denote by $\CD_R^1(\Omega)$ the closure of $ \CD_R^\infty(\Omega)$ in $H^1(\Omega) \times L^2(\Omega)$. In $\CD_R^1(\Omega)$ we consider the energy norm
\begin{equation} \label{eq:energy-norm}
\|(u_0,u_1)\|^2_{E} = \frac{1}{2}\int_{\Omega} |\grad_\g u_0|^2  + c^{-2}q|u_0|^2  + c^{-2}|u_1|^2\dVol.
\end{equation}
When $q$ is not identically zero, the energy norm is equivalent to the $H^1(\Omega)\times L^2(\Omega)$ norm in $\CD^1_R(\Omega)$.
%In other words, there exist a constant $A >0$, such that 
%\begin{equation} \label{eq:norm-equivalence}
%A^{-1}(\|u_0\|^2_{H^1(\Omega)} +  \| u_1\|^2_{L^2(\Omega)})  \leq \|%(u_0,u_1)\|^2_{R} \leq  A(\|u_0\|^2_{H^1(\Omega)} +  \| u_1\|%^2_{L^2(\Omega)}).
%\end{equation}
%for $(u_0, u_1) \in \CD^1$.
However, when $q$ is identically zero,  (\ref{eq:energy-norm}) is not a norm in $\CD_R^1(\Omega)$ since constants are non-trivial solutions of $\Aop u = 0$, $\partial_\nu u + \lambda \partial_t u =0$, and hence $\|(\text{const.},0)\|_{E} =0$. In such case,  we consider the Cauchy data set $\CD_R^{1,0}(\Omega)$, given by the $(u_0,u_1) \in \CD^1(\Omega)$ such that 
\begin{equation} \label{eq:equivalence-condtion}
\int_\Omega u_1(x) dx + \int_{\partial \Omega} \lambda(x)u_0(x) \dS =0
\end{equation}
where $\dS$ is the induced Riemannian metric on $\partial \Omega$ with the induced orientation. It follows easily that (\ref{eq:energy-norm}) is also a norm on $\CD_R^{1,0}(\Omega)$. To justify this condition (\ref{eq:equivalence-condtion}), notice 
\begin{equation} \label{eq:consevation-law}
\int_{\Omega} u_t(t,x) dx + \int_{\partial\Omega}  \lambda(x)u(t,x) dS  = \text{independent of time},
\end{equation}
hence by replacing $u$ in (\ref{Eqn:001f-Cauchy})-(\ref{Eqn:003f-Cauchy}), with 
\[
u - \frac{\int_{\Omega} u_t \, dx + \int_{\partial \Omega} \lambda u \,dS}{\int_{\partial \Omega} \lambda \, dS}   
\]
we reduce the non-uniqueness problem in (\ref{Eqn:001f-Cauchy})-(\ref{Eqn:003f-Cauchy}) to the study of solution for which 
\[
\int_\Omega u_t(t,x) dx + \int_{\partial \Omega} \lambda(x)u(t,x) \dS =0
\]
In summary, we consider  the problem (\ref{Eqn:001f-Cauchy})-(\ref{Eqn:003f-Cauchy}) with Cauchy data set  $\CD_R^1(\Omega)$ when $q$ is not identically zero and $\CD_R^{1,0}(\Omega)$ otherwise; to simplify notation we denote this set by \nolinebreak $\CD_R(\Omega)$. Recall that in $\CD_R(\Omega)$ we always consider the norm $\| \cdot \|_{E}$ given by (\ref{eq:energy-norm}), whether or not $q$ is equals to zero.

For the problem of thermoacoustic and photoacoustic tomography we deal only with the case $u_1 =0$ in (\ref{Eqn:001f-Cauchy})-(\ref{Eqn:003f-Cauchy}). For this reason we denote by $H_R(\Omega)$ the projection of $\CD_R(\Omega)$ onto its first component. The norm induced in $H_R(\Omega)$ is the Dirichlet norm
\begin{equation}\label{eq:dirichlet-norm}
\| u_0 \|^2_{H_R} =  \frac{1}{2}\int_{\Omega} |\grad_\g u_0|^2  + c^{-2}q|u_0|^2  \dVol.
\end{equation}
Clearly $H_R(\Omega)$ is a subspace of $\CD_R(\Omega)$, moreover $\CD_R(\Omega) = H_R(\Omega)\oplus H^0(\Omega)$. The later follows from the fact that the energy of $u$ at time $t$ is given by $E(t,u) =  \| u \|^2_{H_R} + \| u_t\|_{L^2} $.

%In a similar way, but related to Dirichlet boundary conditions, we define $\CD^\infty_D(\Omega)$ and $\CD^1_D(\Omega)$. In this case, $\| \cdot \|_{E}$ is always a norm of $\CD^1_D(\Omega)$ independently of $q$ being zero or not. To be consistent with the notation for the Cauchy data with Robin boundary conditions, we denote by $\CD_D(\Omega)$ the Cauchy data with Dirichlet boundary conditions $\CD^1(\Omega)$ (with no superscript 1).

We now define the evolution operators related to Robin and Neumann boundary conditions. Let $t\geq 0$. First consider the case when $\{ \lambda >0\}$ is not empty. In this case, we denote by $S_R(t): \CD_R(\Omega) \to \CD_R(\Omega)$ the evolution operator of the wave equation related to Robin boundary conditions. The operator $S_R(t)$ takes  Cauchy data at time $0$, to Cauchy data at time $t$. Under the geometric control condition \ref{Assump.001} with exact controllability time $\tau$, this evolution operator $S_R(t)$ defines a contraction semigroup in $\left(\CD_R(\Omega),\|\cdot\|_{E}\right)$, see Theorems 5.5 in \cite{Bar-Leb-Rau-1992} for the case $q$ not zero and Theorem 5.6 in \cite{Bar-Leb-Rau-1992} for the case $q$ equals to zero. Notice that since
\[
\partial_t E(t,u) = - \int_{\partial \Omega} \lambda(x) |u_t|^2 dS,
\]
Robin boundary conditions are energy decreasing only if $\{ \lambda >0 \}$ is non empty. For the case when $\lambda \equiv 0$, denote by $S_N(t): \CD_N(\Omega) \to \CD_N(\Omega)$ the evolution operator of the wave equation related to Neumann boundary conditions. Since Neumann boundary conditions are energy preserving then $\|S_N(t)\| = 1$ as an isomorphism in $\left(\CD_N(\Omega),\|\cdot\|_{E}\right)$.

We now present the Neumann series reconstruction for the case of full data. 

\subsection{Reconstruction Algorithm}

Let $u$ be a solution of the initial value problem (\ref{Eqn:001f-Cauchy})-(\ref{Eqn:003f-Cauchy}). Let $\tau>0$ be the exact controllability time, we define the ``back projection operator'' $A: H^1((0,\tau)\times \partial \Omega) \to H^1(\Omega)$ that acts on the data $\Lambda u_0 = u|_{(0,\tau)\times \partial \Omega}$ by 
\[
(A \Lambda u_0)(x) = v(0,x)
\]
where $v$ solves 
\numparts 
\begin{eqnarray}
\partial_{t}^{2} {v} - \Aop v = 0 \qquad && \mbox{in} \quad (0,\tau) \times \Omega  \\
v = 0 \quad \mbox{and} \quad \partial_{t} v = 0 \qquad && \mbox{on} \quad \{ t = \tau \} \times \Omega  \\
\partial_{\nu} v  = - \lambda \partial_{t} \Lambda u_0   \qquad && \mbox{on} \quad (0,\tau) \times  \Gamma \\
\partial_{\nu} v  = 0   \qquad && \mbox{on} \quad (0,\tau) \times \partial \Omega \setminus \Gamma 
\end{eqnarray}
\endnumparts
We think of $A\Lambda u_0$ as the first approximation on the reconstruction of $u_0$. Notice that we only need the data $\Lambda u_0$ on $\Gamma$. The  ``error operator''  $K$ is given by 
\[
K := \text{Id} - A\Lambda
\]
Notice that if we let $w = u - v$, then $(K u_0)(x) = w(0,x)$; moreover $w$ solves
\numparts 
\begin{eqnarray}
\partial_{t}^{2} {w} - \Aop w = 0 \qquad && \mbox{in} \quad (0,\tau) \times \Omega  \\
w = u(\tau, \cdot)  \quad \mbox{and} \quad \partial_{t}  w =  \partial_t u(\tau,\cdot) \qquad && \mbox{on} \quad \{ t = \tau \} \times \Omega  \\ 
\partial_{\nu} w =  0 \qquad && \mbox{on} \quad (0,\tau) \times \partial \Omega 
\end{eqnarray}
\endnumparts 
Then we can write $K: H_R(\Omega) \to H_R(\Omega)$ as a composition of operators
\[
K = \pi_1 S_N(-\tau) S_R(\tau) \pi^*_1 
\]
where $\pi^*(u_0) = (u_0,0)$ and $\pi_1(u_0,u_1) = u_0$.

The important observation is that, if $\Gamma = \{\lambda >0\}$ satisfies the GCC \ref{Assump.001} with exact controllability time $\tau$, $K$ is the composition of unitary operators and $S_R(\tau)$ which is a contraction, see Theorems 5.5 and 5.6 in  \cite{Bar-Leb-Rau-1992}. Hence $K$ itself is a contraction, i.e., $\|K\|_{H_R \to H_R} < 1$. This implies in particular that we can reconstruct $u_0$ from the data $\Lambda u_0|_{\Gamma}$ via a Neumann series expansion. We state this observation as a theorem. 

\begin{theorem}[\textbf{Reconstruction Method 2}] 
\label{Thm.Reconstruction2}
Let $\Gamma = \{\lambda >0\}$ be non empty and satisfying the GCC \ref{Assump.001} with exact controllability time $\tau$. Then $A\Lambda = \text{Id} - K$, where $K$ is a contraction, $\|K\|_{H_R \to H_R} < 1$. In particular, $\text{Id} - K$ is invertible on $H_R(\Omega)$ and the inverse thermoacoustic problem has an explicit Neumann series reconstruction given by 
\begin{equation}\label{eq:Neumann-series}
u_0 = \sum_{n=0}^\infty K^nA\Lambda u_0.
\end{equation}
\end{theorem}

%%%%%%%%%%%%%%%%%%%%%%%%%%%%%%%%%%%%%%%%%%%%%%%%%%%%%%%%%%%%%%%%%
%%%%%%% NEW SECTION %%%%%%%%%%%%%%%%%%%%%%%%%%%%%%%%%%%%%%%%%%%%%
%%%%%%%%%%%%%%%%%%%%%%%%%%%%%%%%%%%%%%%%%%%%%%%%%%%%%%%%%%%%%%%%%

\section{Conclusion} \label{Section:Discussion}

We have proposed two reconstruction algorithms for the thermoacoustic tomography problem in an enclosure $\Omega$ with observation on a portion $\Gamma$ of the enclosing boundary $\partial \Omega$. In both cases, the reflecting portion $\partial \Omega \setminus \Gamma$ of the boundary allows the energy of the unknown initial condition to eventually reach $\Gamma$ where measurements take place. Our work provides precise conditions under which such a scenario yields the solvability of this inverse problem, and explains from a mathematical viewpoint how the use of acoustic mirrors may render a good reconstruction for limited-view detection.      

The first algorithm (Theorem \ref{Thm.Reconstruction}), based on the conjugate gradient method, is naturally suited for \textit{full enclosure} where the acoustic energy may not be allowed to escape from the domain of interest. The second algorithm (Theorem \ref{Thm.Reconstruction2}), based on a Neumann series expansion, is computationally more attractive and naturally suited for \textit{partial enclosure} where the energy is allowed to leave through the observation part of the boundary. The authors are in the process of developing numerical implementations to compare the proposed reconstruction methods with each other and with experimental data. As soon as meaningful results are obtained from these efforts, they will be reported in a forthcoming publication.

%% %%%%%%%%%%%%%%%%%%%%%%%%%%%%%%%%%%%%%%%%%%%%%%%%%%%%%%%%%%%%%%%%%%%%%%%%%%%%%%%

\section*{Acknowledgments} \label{Sec:Acknowledgements}
The authors would like to thank Plamen Stefanov and Yang Yang for fruitful conversations and for providing a draft of their paper \cite{stefanov-yang-2014}. We are also grateful for recommendations made by Gunther Uhlmann.

%% %%%%%%%%%%%%%%%%%%%%%%%%%%%%%%%%%%%%%%%%%%%%%%%%%%%%%%%%%%%%%%%%%%%%%%%%%%%%%%%

\section*{References}

\bibliographystyle{plain}
\bibliography{Biblio}

\end{document}